# Bornes effectives pour la torsion des courbes elliptiques sur les corps de nombres

Pierre Parent

8 novembre 1995

**Résumé**

On se propose de donner une forme effective au théorème de Mazur-Kamienny-Merel sur la torsion des courbes elliptiques sur les corps de nombres. Pour ce faire, on prouve ici qu'est satisfait un critère de Kamienny qu'on explicite en niveau non premier.

## 1 Introduction

*La "conjecture de borne uniforme pour les courbes elliptiques", affirmant qu'il existe pour tout entier $d$ un entier $B(d)$ tel que, pour tout corps de nombres $K$ de degré $d$ sur $\mathbb{Q}$ et pour toute courbe elliptique $E$ sur $K$, la partie de torsion $E(K)_{tors}$ du groupe de Mordell-Weil $E(K)$ est de cardinal majoré par $B(d)$, a été démontrée dans le cas général en février 1994 par Loïc Merel. En fait, Merel (et Oesterlé) montrent que, si $P$ est un point d'ordre $p$ premier de $E(K)$, on a $p \leq (1 + 3^{d/2})^2$. Des travaux de Faltings et Frey permettent alors de conclure à l'existence des bornes $B(d)$, mais pas de manière effective : en effet, si on a bien majoré les nombres premiers pouvant diviser les groupes $E(K)_{tors}$, on ne sait pas en pratique quelles puissances de ces nombres premiers peuvent intervenir dans ces groupes. Le but de cet article est de donner une preuve du fait qu'est satisfait, en niveau $p^n$ "explicitement assez grand" par rapport à $d$, un "critère de Kamienny" (Proposition 4), qui impliquera à son tour une telle borne (voir le Corollaire 6). La rédaction complète de la preuve de cette implication, (seulement esquissée ici, et essentiellement identique à celle du cas où le niveau est premier), fera l'objet d'un prochain article.*

## 2 Schéma de la preuve

Soit $E$ une courbe elliptique définie sur un corps de nombres $K$ de degré $d$ sur $\mathbb{Q}$, et possédant un point $K$-rationnel $P$, tel que le cardinal du sous-groupe cyclique de $E(K)$ engendré par $P$ soit une puissance d'un nombre premier : $|\langle P \rangle| = p^n$. On cherche à majorer $p^n$ en fonction de $d$. Soit $l$ un nombre premier différent de 2 et de $p$ (on prendra en pratique le plus petit possible, *i.e.* $l = 5$ si $p = 3$, $l = 3$ dans tous les autres cas). Soit $\mathcal{L}$ une place de $\mathcal{O}_K$ au-dessus de $l$. Examinons la fibre en $\mathcal{L}$ du modèle de Néron de la courbe elliptique (qu'on note encore $E$) sur $\mathcal{O}_K$ (l'anneau des entiers du corps $K$).

**Proposition 1** *Si en $\mathcal{L}$ au-dessus de $l$, on a l'un des trois cas:*

1. *$E$ a bonne réduction;*

2. *$E$ a réduction additive;*

3. *$E$ a réduction multiplicative et $\langle \widetilde{P} \rangle$ (où $\widetilde{P}$ est la réduction de $P$ mod $\mathcal{L}$) est inclus dans sa composante neutre,*

*alors $|\langle P \rangle| \leq (1 + l^{d/2})^2$.*

Supposons donc qu'en tout $\mathcal{L}$ au-dessus de $l$, $E$ ait réduction multiplicative et que $\widetilde{P}$ ne soit pas trivial dans le groupe des composantes de cette réduction $\widetilde{E}$ de $E$. On aurait aimé pour la suite que $\widetilde{P}$ soit d'ordre $p^n$ dans le groupe des composantes de $\widetilde{E}$: d'où l'idée d'examiner le quotient de $\widetilde{E}$ par le plus gros sous-groupe de $\langle \widetilde{P} \rangle$ inclus dans la composante neutre de $\widetilde{E}$. Soit donc $n_\mathcal{L}$ la différence entre $n$ et l'ordre de ce plus gros groupe: $n_\mathcal{L}$ est le plus petit entier tel que $p^{n_\mathcal{L}}.\widetilde{P}$ tombe dans dans la composante neutre de $\widetilde{E}$. Soit aussi $n'$ le plus petit des $n_\mathcal{L}$, pour $\mathcal{L}$ parcourant l'ensemble des places de $\mathcal{O}_K$ au-dessus de $l$. Oesterlé démontre:

**Lemme 2** *Soit $k$ élément de $\mathbb{N}$; si en la place $\mathcal{L}$ de $\mathcal{O}_K$, $p^{k-1}P$ ne se réduit pas dans la composante neutre de $\widetilde{E}$, alors la réduction de $\widetilde{P}$ dans $\widetilde{E}/\langle p^k.\widetilde{P}\rangle$ est d'ordre exactement $p^k$ dans son groupe des composantes.*

L'interprétation modulaire de ce lemme est donc la suivante: le point $K-$rationnel $j$ de la courbe modulaire $X_0(p^n)$ que définit le couple $(E/\langle p^{n'}.P \rangle, \langle P \rangle)$ se réduit en la pointe 0 modulo toute place $\mathcal{L}$; et l'image $j'$ de ce point par l'involution d'Atkin-Lehner, en la pointe infinie. Voyons comment on peut alors utiliser les arguments de Kamienny ([Kamienny 92b], [Edixhoven 94]) pour démontrer son critère.

Si les $\sigma_i$, $1 \leq i \leq d$ sont les plongements de $K$ dans $\overline{\mathbb{Q}}$, $j'^{(d)} = (\sigma_1(j'), \sigma_2(j'), ..., \sigma_d(j'))$ définit un point $\mathbb{Q}$-rationnel du produit symétrique $d$-ième: $X_0(p^n)^{(d)}$, de $X_0(p^n)$.

On utilise comme Merel le quotient d'enroulement ([Merel 95]). On considère les premiers groupes d'homologie singulière absolue: $H_1(X_0(p^n);\mathbb{Z})$ et relative aux pointes: $H_1(X_0(p^n), \text{pointes};\mathbb{Z})$, de $X_0(p^n)$. Si $a$ et $b$ sont deux éléments de $\mathbb{P}^1(\mathbb{Q})$, le *symbole modulaire* $\{a,b\}$ est l'élément de $H_1(X_0(p^n), \text{pointes};\mathbb{Z})$ défini par l'image de n'importe quel chemin continu reliant $a$ à $b$ sur le demi-plan de Poincaré auquel on a ajouté l'ensemble $\mathbb{P}^1(\mathbb{Q})$ de ses pointes (dont l'ensemble des points complexes de la courbe modulaire, vue comme surface de Riemann, est un quotient). On a un isomorphisme d'espaces vectoriels réels:

$$\begin{cases} H_1(X_0(p^n);\mathbb{Z}) \otimes \mathbb{R} \to \text{Hom}_\mathbb{C}\left(H^0(X_0(p^n);\Omega^1),\mathbb{C}\right) \\ \gamma \otimes 1 \mapsto \left(\omega \mapsto \int_\gamma \omega\right) \end{cases}$$

Selon un théorème de Manin, l'image réciproque de la forme linéaire $\omega \mapsto -\int_{\{0,\infty\}} \omega$ dans $H_1(X_0(p^n);\mathbb{R})$ est en réalité dans $H_1(X_0(p^n);\mathbb{Q})$. Notons-la $e$ (comme d'habitude). On note (toujours comme d'habitude) $\mathbb{T}$ l'algèbre engendrée par les opérateurs de Hecke $T_i$, $i \in \mathbb{N}^*$, agissant entre autres sur $H_1(X_0(p^n); \mathbb{Q})$ et sur la jacobienne $J_0(p^n)$ de la courbe modulaire. Soit $\mathcal{A}_e$ l'idéal annulateur dans $\mathbb{T}$ de $e$ (*idéal d'enroulement*); on définit alors le *quotient d'enroulement* $J_0^e$ comme la variété abélienne quotient $J_0/\mathcal{A}_e J_0$. De façon similaire à Merel dans [Merel 95], on montre la:

**Proposition 3** $J_0^e(\mathbb{Q})$ *est fini.*

Soit maintenant l'application naturelle $f_d: X_0(p^n)_{\text{lisse}}^{(d)} \to J_0^e$, qu'on a normalisée par $\infty^{(d)} \mapsto 0$. On montre "par les arguments standards" que le fait que $f_d$ soit une immersion formelle en $\infty_l^{(d)}$ contredit l'existence de notre point $j'^{(d)}$. Or on a le "critère de Kamienny":

**Proposition 4** *On a équivalence entre:*
*1)* $f_d$ *est une immersion formelle en* $\infty_l^{(d)}$, *et*
*2)* $T_1 e, ..., T_d e$ *sont* $\mathbb{F}_l$-*linéairement indépendants dans* $\mathbb{T}e/l\mathbb{T}e$.
*De plus, ces deux conditions sont satisfaites si l'est:*
*3)* $T_1\{0, \infty\}, ..., T_{d.s}\{0, \infty\}$ *sont* $\mathbb{F}_l$-*linéairement indépendants dans* $H_1(X_0(p^n), \text{pointes}; \mathbb{Z})$
*(ici et pour la suite, s désigne le plus petit nombre premier différent de p).*

(Le fait que la dernière condition implique les précédentes est une remarque d'Oesterlé). Il suffit donc maintenant de prouver:

**Proposition 5** *Soit* $C = (4096.\pi^2)/(2\sqrt{2} - 1)$. *Si* $p^n > C^2.(sd)^6$, *alors les* $T_i\{0,\infty\}$, $1 \leq i \leq sd$ *sont* $\mathbb{F}_m$-*linéairement indépendants (dans* $H_1(X_0(p^n), \text{pointes}; \mathbb{Z}) \otimes \mathbb{F}_m$) *pour tout corps fini* $\mathbb{F}_m$ *avec* $m$ *premier.*

De cette proposition découle en effet:

**Corollaire 6** *Soit* $E$ *une courbe elliptique sur un corps* $K$ *de degré* $d$ *sur* $\mathbb{Q}$. *Si* $E(K)$ *possède un point* $P$ *d'ordre une puissance* $p^n$ *d'un nombre premier* $p$, *on a:*

1. $p^n \leq C^2.(3^d + 1).(2d)^6$, *si* $p$ *est différent de* 2 *et* 3;

2. *Si* $p = 3$, $p^n \leq C^2.(5^d + 1).(2d)^6$;

3. *et pour* $p = 2$, $2^n \leq C^2.(1 + 3^d).(3d)^6$,

*où* $C$ *est le réel défini dans la proposition 5.*

**Remarques** *On sait qu'on a* $E(K)_{tors} \cong \mathbb{Z}/n_1\mathbb{Z} \times \mathbb{Z}/n_2\mathbb{Z}$. *On borne donc ainsi le cardinal de tous les p-groupes de* $E(K)_{tors}$, *et avec la borne de Merel-Oesterlé pour les nombres premiers p pouvant intervenir on obtient une borne effective "globale" pour l'ensemble* $\{card(E(K)_{tors}) \mid [K:\mathbb{Q}] = d \text{ et } E \text{ est une courbe elliptique sur } K\}$.
*Une telle borne globale semble de toute façon assez facile à améliorer en reprenant les arguments de ce papier directement en niveau* $N$ *entier quelconque.*

La preuve dans son ensemble ayant été esquissée, donnons maintenant la démonstration complète de la dernière proposition, qui permet d'appliquer le critère de Kamienny.

## 3    Preuve de la Proposition 5

On identifie $\Gamma_0(p^n)\backslash SL_2(\mathbb{Z})$ à $\mathbb{P}^1(\mathbb{Z}/p^n\mathbb{Z})$ avec :

$$\Gamma_0(p^n)\begin{pmatrix} a & b \\ c & d \end{pmatrix} \mapsto (\overline{c},\overline{d}) = (c\bmod(p^n), d\bmod(p^n))$$

L'application de $\Gamma_0(p^n)\backslash SL_2(\mathbb{Z})$ vers $H_1(X_0(p^n),\text{ pointes};\mathbb{Z}) : g \mapsto \{g\cdot 0, g\cdot\infty\}$ s'identifie alors a une application de $\mathbb{P}^1(\mathbb{Z}/p^n\mathbb{Z})$ vers la même chose, qu'on note $\xi$ :

$$\xi(\overline{w},\overline{t}) = \left\{ \overline{\begin{pmatrix} a & b \\ w & t \end{pmatrix}\cdot 0}, \overline{\begin{pmatrix} a & b \\ w & t \end{pmatrix}\cdot\infty}\right\} = \left\{\frac{b}{t},\frac{a}{w}\right\},$$

avec $w,t$, relèvements dans $\mathbb{Z}$ de $\overline{w}$ et $\overline{t}\in\mathbb{Z}/p^n\mathbb{Z}$ et $a,b\in\mathbb{Z}$ tels que $\begin{pmatrix} a & b \\ w & t \end{pmatrix}$ soit dans $SL_2(\mathbb{Z})$.

On sait de plus qu'on a :

$$T_r\{0,\infty\} = \sum_{\substack{0\leq w < t \\ 0\leq v < u \\ ut-vw=r}} \xi(\overline{w},\overline{t}),$$

où on pose $\xi(\overline{w},\overline{t}) = 0$ si $\text{pgcd}(w,t,p) > 1$.

Soit $\sigma = \overline{\begin{pmatrix} 0 & 1 \\ -1 & 0 \end{pmatrix}}$ et $\tau = \overline{\begin{pmatrix} 0 & -1 \\ 1 & -1 \end{pmatrix}}$. On choisit pour représentants de $\mathbb{P}^1(\mathbb{Z}/p^n\mathbb{Z})$ : $\{(R_1,1),\ R_1$ un système de représentants de $\mathbb{Z}/p^n\mathbb{Z}\}\ \cup\ \{(1,p.R_2),\ R_2$ un système de représentants de $\mathbb{Z}/p^{n-1}\mathbb{Z}\}$. On note $w/t$ au lieu de $(\overline{w},\overline{t})$, souvent. On fait agir $SL_2(\mathbb{Z})$ sur $\mathbb{P}^1(\mathbb{Z}/p^n\mathbb{Z})$ à droite, comme d'habitude. En particulier :

$$(\overline{w},\overline{t})\cdot\sigma = (\overline{w},\overline{t})\overline{\begin{pmatrix} 0 & -1 \\ 1 & 0 \end{pmatrix}} = (\overline{-t},\overline{w}),$$

et de même $(w/t)\cdot\tau = -t/(w+t)$.

On note encore $\Sigma_r = \{(\overline{w},\overline{t}),\ 0\leq w < t\ /\ $il existe $(u,t),\ 0\leq u < v,\ 0\leq (ut-vw)\leq r\}\backslash\{(\overline{1},\overline{r})\}$, et $\mathbb{Z}[\mathbb{P}^1(\mathbb{Z}/p^n\mathbb{Z})]^\sigma$ (respectivement, $\mathbb{Z}[\mathbb{P}^1(\mathbb{Z}/p^n\mathbb{Z})]^\tau$) désigne l'ensemble des éléments de $\mathbb{Z}[\mathbb{P}^1(\mathbb{Z}/p^n\mathbb{Z})]$ stables par l'action de $\sigma$ (respectivement, $\tau$). Le symbole $\sum_\sigma$ désignera une somme à valeurs dans $\mathbb{Z}[\mathbb{P}^1(\mathbb{Z}/p^n\mathbb{Z})]^\sigma$, et de même avec $\tau$.

On a la suite exacte ([Merel 93]) :

$$0\to\mathbb{Z}\stackrel{\phi_1}{\to}\mathbb{Z}[\mathbb{P}^1(\mathbb{Z}/p^n\mathbb{Z})]^\sigma\times\mathbb{Z}[\mathbb{P}^1(\mathbb{Z}/p^n\mathbb{Z})]^\tau\stackrel{\phi_2}{\to}$$

$$\stackrel{\phi_2}{\to}\mathbb{Z}[\mathbb{P}^1(\mathbb{Z}/p^n\mathbb{Z})]\stackrel{\phi_3}{\to} H_1(X_0(p^n),\text{ pointes};\mathbb{Z})\to 0,$$

avec $\phi_1 : n \mapsto (n\sum_{\mathbb{P}^1} x, -n\sum_{\mathbb{P}^1} x)$, $\phi_2 : (\Sigma\alpha_x x, \Sigma\beta_x x) \mapsto (\Sigma\alpha_x x + \Sigma\beta_x x)$, et $\phi_3 : \Sigma\lambda_x x \mapsto \Sigma\lambda_x.\xi(x)$.

**Preuve de la proposition** Posons $D = sd$. Supposons que l'on ait une relation de liaison: $\sum_{i=1}^{r} \lambda_i T_i\{0,\infty\} = 0$ dans $H_1(X_0(p^n), \text{pointes}; \mathbb{F}_l)$, avec $r \leq D$. On va montrer que $\lambda_r$ est nul: on aura ainsi le résultat par récurrence immédiate.

Ce qui précède permet d'écrire:

$$\sum_{i=1}^{r} \lambda_i T_i\{0,\infty\} = 0 \iff \lambda_r.(1/r) + \sum_{\Sigma_r} \mu_{(\overline{w},\overline{t})}(w/t) = \sum_\sigma \alpha_{(\overline{w},\overline{t})}(\overline{w},\overline{t}) - \sum_\tau \beta_{(\overline{w},\overline{t})}(\overline{w},\overline{t}).$$

De même que dans [Merel 93] on prouve avec l'aide de Fouvry (voir plus loin) le:

**Lemme 7** *Si $A$ et $B$ sont deux intervalles de $\{1, 2, ..., p^n - 1\}$, et:*

1. $|A|, |B| \geq \sup(11, 2p+1)$,

2. $(|A| - 11).(|B| - 11) > 144 + (C/8).p^{3n/2}$,

*$C$ étant le réel défini dans la Proposition 5, alors il existe $y \in A$ et $z \in B$ tel que $y.z = -1 \mod (p^n)$.*

On considère le graphe de $\overline{\mathbb{P}^1(\mathbb{Z}/p^n\mathbb{Z})}$ dont les arêtes sont l'action de $\sigma$ et $\tau$ (voir **Figure 1**): on a $(\tau\sigma) = \begin{pmatrix} -1 & 0 \\ -1 & -1 \end{pmatrix}$, donc $(\overline{w},\overline{t}) \cdot \tau\sigma = (\overline{w},\overline{t}) + \overline{1}$ (et $(\overline{w},\overline{t}) \cdot \sigma\tau^2 = (\overline{w},\overline{t}) - \overline{1}$).

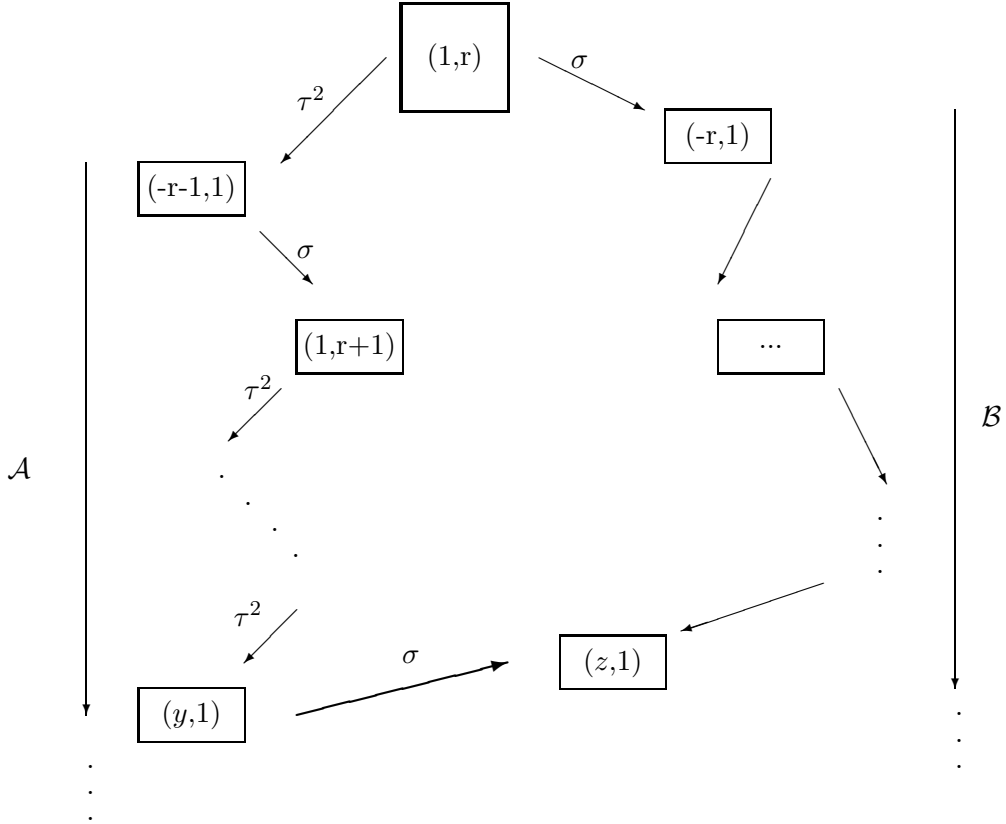

**Figure 1**

On va montrer qu'on a sur ce graphe, "de part et d'autre" de $(\overline{1}, \overline{r})$ (c'est-à-dire, contenant $(\overline{1}, \overline{r}) \cdot \tau^2 = (\overline{-r-1}, \overline{1})$ et $(\overline{1}, \overline{r}) \cdot \sigma = (\overline{-r}, \overline{1})$ respectivement), deux chemins $\mathcal{A}$ et $\mathcal{B}$ ne rencontrant pas d'éléments de $\Sigma_r$, (somme des autres éléments du graphe intervenants dans $\sum \lambda_i T_i\{0, \infty\}$), et qui contiennent des intervalles de cardinal supérieur à $(p^n/D) - D - 2$ et $(p^n/D^2) - 2$. Alors par le lemme 7, pour
$$p^{2n}/D^3 > C.p^{3n/2} ,$$
i.e.
$$p^n > C^2.D^6 ,$$
on aura dans $\mathcal{A}$ et dans $\mathcal{B}$ respectivement des éléments $y$ et $z$ tel que $y \cdot \sigma = \frac{-1}{y} = z$, et $y \cdot \sigma$ sera un élément de $\mathcal{A}$, puisque les deux seront de forme :

$$... \to (\overline{a}, \overline{1}) \xrightarrow{\tau} (\overline{-1}, \overline{a+1}) \xrightarrow{\sigma} (\overline{a+1}, \overline{1}) \to ...$$
$$\longrightarrow \quad +\overline{1} = \tau\sigma \quad \longrightarrow$$

Les deux chemins se "rencontrent", donc. Or on a, pour tout élément $x$ du graphe, $\mu_x = \alpha_x - \beta_x$ par ce qui précède. Puisque les deux chemins ne rencontrent pas $\Sigma_r$, si $x$ en est, on a $\mu_x = 0$ donc $\alpha_x = \beta_x$. De plus, on parcourt ces chemins en faisant agir $\sigma$ ou $\tau$; donc si $x' = x \cdot \sigma$, $\alpha_{x'} = \alpha_{x \cdot \sigma} = \alpha_x = \beta_x = \beta_{x'}$, de même avec $\tau$ - donc $\alpha_x = \beta_x \equiv \alpha_{-r}$ sur les deux chemins. Mais

$$\lambda_{\overline{\frac{1}{r}}} = \mu_{\frac{1}{r}} = \alpha_{\frac{1}{r}} - \beta_{\frac{1}{r}} = \alpha_{\frac{1}{r} \cdot \sigma} - \beta_{\frac{1}{r} \cdot \tau^2} = 0.$$

On montre donc l'existence de ces chemins. (Dans tous les calculs qui suivent, on confond l'écriture d'un entier $w$ et de sa réduction $\overline{w}$, pour alléger les notations) .

<u>Premier chemin :</u> $\mathcal{A}$ partant de $\frac{1}{r} \tau^2 = -r - 1$

1) Si $\frac{w}{t} - (-r - 1) = \overline{a}$, avec $a$ choisi dans $\{-p^n + 1, ..., -1, 0\}$, et $\overline{a}$ : classe de $a \bmod p^n$ ( $\iff w + t(r + 1) = at + b.p^n$, $b \in \mathbb{Z}$).

<u>Si $b = 0$</u> $t(r + 1) - at = -w$ : incompatibilité de signes.

<u>Si $b \neq 0$</u> $|a| = \frac{1}{t}|b.p^n - w - t(r+1)| \geq \frac{(p^n - D(D+1) - D)}{D} \geq \frac{p^n}{D} - D - 2$ (on a en effet : $det \begin{pmatrix} u & v \\ w & t \end{pmatrix} = k \leq r$, $u > v \geq 0$, $t > w \geq 0$, donc : $k = ut - vw \geq ut - (u-1)(t-1)$ et $u + t - 1 \leq r$, $t \leq r - u + 1 \leq r \leq D$).

2) Si $(\frac{w}{t})\sigma - (-r - 1) = \frac{-t}{w} + r + 1 = \overline{a}$, $a \in \{-p^n + 1, ..., -1, 0\}$ ; $-t + w(r + 1) = aw + bp^n$ ;

<u>Si $b = 0$</u> $w(r + 1 - a) = t$ ; mais $(r + 1 - a) \geq r + 1$, et $0 \leq w < t \leq r$ : contradiction.

<u>Si $b \neq 0$,</u> $|a| \geq (\frac{p^n - D(D+1) + D}{D}) \geq \frac{p^n}{D} - D$.

On peut donc "reculer" (... $\alpha \xrightarrow{\sigma} . \xrightarrow{\tau^2} \alpha - 1 \xrightarrow{\sigma} . \xrightarrow{\tau^2} \alpha - 2$ ...) à partir de $(-r - 1)$, et décrire ainsi un chemin contenant un intervalle de cardinal supérieur à $p^n/D - D - 2$.

<u>Second chemin</u> On doit là distinguer deux cas, selon que $p$ divise ou non $r$ (voir **Figure 2**).

**a)** Si $p$ ne divise pas $r$, chemin $\mathcal{B}$ : on part de $(\frac{1}{r})$ lui même, on recule de même :

1) $\frac{w}{t} - \frac{1}{r} = \overline{a}$, $-p^n < a \leq 0 \iff wr - t = art + b.p^n$ ;

<u>$b = 0$</u> : $t = r(w - at) \Rightarrow a = 0$, $w = 1$, $t = r$ : c'est $\frac{1}{r}$ lui-même.

<u>$b \neq 0$</u> : $|a| \geq \frac{p^n - D^2}{D^2}$ de même que plus haut.

2) $-\frac{t}{w} - \frac{1}{r} = \overline{a}$ : $-rt - w = awr + b.p^n$ ;

$\underline{b = 0} \Rightarrow r|w$, impossible (car $w \leq r - 1$, et : $w = 0 \Rightarrow t = 0$, or $t > 0$).

$\underline{b \neq 0} \Rightarrow |a| \geq \frac{p^n}{D^2} - 2$

**b)** Si $p$ divise $r$, on a alors que le chemin $\mathcal{B}$ précédent : $\frac{1}{r} \xrightarrow{\sigma\tau^2} \frac{1-r}{r} \xrightarrow{\sigma\tau^2} ...$
est bien de longueur supérieure à $(\frac{p^n}{D^2})$ ; mais il ne contient pas cette fois d'intervalle, puisque $r$ n'est pas inversible modulo $p^n$, donc l'action de $\sigma\tau^2$ ne correspond plus à l'addition de $(-1)$ (les $k/pl$ ne sont plus relevables en éléments de $\mathbb{Z}/p^n\mathbb{Z}$). Cependant, le calcul précédent montre que $\frac{r}{r-1} = \frac{1}{r} \cdot \sigma\tau^2\sigma$ ; et $(r-1)$ est, cette fois, inversible modulo $p^n$, donc en "avançant" à partir de cet élément et à l'aide de $(\tau\sigma)$, on aura bien un intervalle ; on note $\mathcal{B}'$ ce chemin. On minore encore une fois sa longueur :

1) $\frac{w}{t} - \frac{r}{r-1} = \overline{a} \iff w(r-1) - rt = at(r-1) + b.p^n$ (on écrit cette fois cela avec $0 \leq a \leq p^n - 1$, et $b \in \mathbb{Z}$)

$\underline{\text{Si } b = 0}$   $t(r + a(r-1)) = w(r-1)$ ; mais $w < t$, et $a(r-1) + r > r - 1$, contradiction.

$\underline{\text{Si } b \neq 0}$   $|a| \geq \frac{p^n}{D^2} - 1$

2) $\frac{w}{t} \cdot \sigma - \frac{r}{r-1} = -\frac{t}{w} - \frac{r}{r-1} = \overline{a} \iff -t(r-1) - rw = aw(r-1) + b.p^n$

$\underline{b = 0}$   $t(r-1) = -w(a(r-1) + r)$, contradiction de signes.

$\underline{b \neq 0}$   $|a| \geq \frac{p^n}{D^2} - 2$

Dans chaque cas, on obtient bien que ce second chemin contient un intervalle de cardinal supérieur ou égal à $(p^n/D^2) - 2$.□

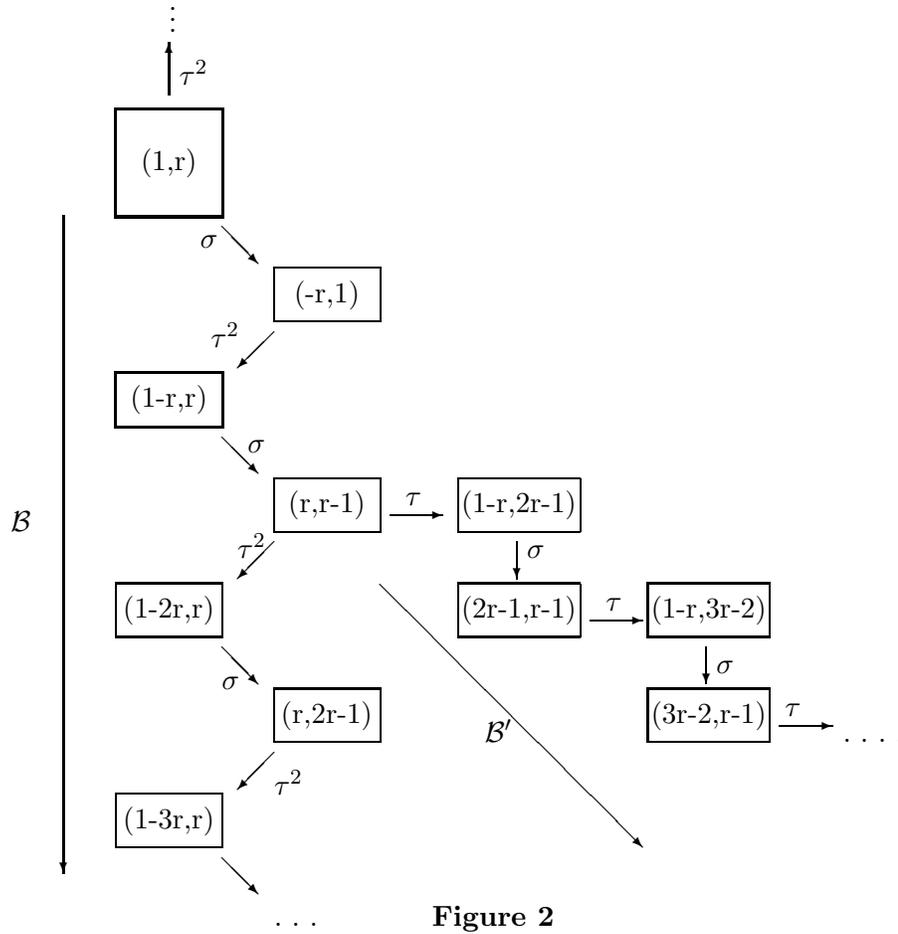

**Figure 2**

# 4 Annexe : preuve du lemme de théorie analytique des nombres

On rappelle son énoncé :
**Lemme 7** *Si $A$ et $B$ sont deux intervalles de $\mathbb{Z}/p^n\mathbb{Z}$ ($p$ un nombre premier), et :*

1. $|A|, |B| \geq \max(11, 2p+1)$,

2. $(|A|-11).(|B|-11) > 144 + (C/8).p^{3n/2}$,

*où $C = (4096.\pi^2)/(2\sqrt{2}-1)$, alors il existe $y \in A$ et $z \in B$ tel que $y.z = -1 \mod (p^n)$.*

**Preuve** On suppose (quitte à lui enlever un élément à l'une de ses extrémités) que le cardinal de $A$ est pair, et on écrit $A = \{a+1, a+2, ..., a+2K\}$ (respectivement $B = \{b+1, b+2, ..., b+2K'\}$ ; on a donc en réalité $K = \frac{|A|}{2}$ ou $\frac{|A|-1}{2}$, et de même avec $K'$). Soit $\Psi_A$ la fonction de $\mathbb{Z}/p^n\mathbb{Z}$ dans $\mathbb{N}$, affine par morceaux et à support dans $A$, définie par :

$$\begin{cases} \Psi_A(a+k) = k/K \text{ si } 1 \leq k \leq K, \\ \Psi_A(a+2K+1-k) = k/K \text{ si } 1 \leq k \leq K, \\ \Psi_A(x) = 0 \text{ si } x \text{ n'appartient pas à } A. \end{cases}$$

On définit de même $\Psi_B$. Il suffit de maintenant prouver que l'expression :

$$\Lambda = \sum_{r \text{ inversible} \mod (q)} \Psi_A(r).\Psi_B(\overline{r})$$

est strictement positive (on a écrit $q = p^n$ et $\overline{r} = -1/r \mod (q)$). En prenant les transformées de Fourier des fonctions, on a :

$$\Lambda = \sum_{r \text{ inversible} \mod q} \left[\frac{1}{q}\sum_{h=0}^{q-1}(\sum_{m=0}^{q-1}\Psi_A(m)\,e^{2i\pi\frac{mh}{q}})e^{-2i\pi\frac{h}{q}r}\right]\left[\frac{1}{q}\sum_{h'=0}^{q-1}(\sum_{m'=0}^{q-1}\Psi_B(m')e^{2i\pi\frac{m'h'}{q}})e^{-2i\pi\frac{h'}{q}\overline{r}}\right]$$

$$\Lambda = \frac{1}{q^2}\sum_{h=0}^{q-1}\sum_{h'=0}^{q-1}\left(\sum_{m\in A}\Psi_A(m)e^{\frac{2i\pi hm}{q}}\right)\left(\sum_{m'\in B}\Psi_B(m')e^{\frac{2i\pi h'm'}{q}}\right)S(-h,-h';q)$$

où :

$$S(-h,-h';q) = \sum_{\substack{r=0 \\ p \text{ ne divise pas } r}}^{q-1} e^{\frac{2i\pi}{q}(-h.r-h'.\overline{r})}$$

est une somme de Kloosterman, qui vérifie :
$S(0,0;q) = \phi(q) = p^n - p^{n-1}$ ;
$S(-h \neq 0, 0; q) = \sum_{r=0}^{p^n-1} e^{2i\pi(-hr)/p^n} - \sum_{r=0}^{p^{n-1}-1} e^{2i\pi(-hr)/p^{n-1}}$
$= \begin{cases} 0 \text{ si } h \neq 0 \mod p^{n-1} \\ -p^{n-1} \text{ si } h = 0 \mod p^{n-1} \end{cases}$, et en général :
$|S(h \neq 0, h' \neq 0; q)| \leq 2\sqrt{2}(-h, -h', q)^{\frac{1}{2}}.\sqrt{q}$

où $(h, h', q)$ désigne le plus grand diviseur commun à $h, h'$ et $q$ (voir [Hooley], [Salié]).

On isole dans $\Lambda$ différents termes:
Terme principal, en $(h = 0, \ h' = 0)$:

$$\Lambda_{0,0} = \frac{1}{q^2} S(0,0;q) \left( \sum_m \Psi_A(m) \right) \left( \sum_{m'} \Psi_B(m') \right) = \frac{p^n - p^{n-1}}{p^{2n}} \left( \frac{2}{K} \sum_0^K k \right) \left( \frac{2}{K'} \sum_0^{K'} k' \right)$$

$$= \frac{(1 - \frac{1}{p})}{p^n}(K+1)(K'+1)$$

et donc
$$\Lambda_{0,0} \geq \frac{1}{2q}(K+1)(K'+1).$$

Somme des termes en $(h \neq 0, h' = 0)$

$$\begin{aligned}
\Lambda_{0,} &= \sum_{h=1}^{q-1} \frac{1}{q^2} \sum_m \Psi_A(m) e^{\frac{2i\pi hm}{q}} \left( \sum_{m'} \Psi_B(m') \right) S(-h, 0 \ ; q) \\
&= \sum_{l=1}^{p-1} \frac{1}{q^2} \sum_m \Psi_A(m) e^{\frac{2i\pi lm}{p}} \left( \sum_{m'} \Psi_B(m') \right) (-p^{n-1}) \\
&= (-1/p^{n+1})(K'+1) \sum_m \Psi_A(m) \sum_{l=1}^{p-1} e^{2i\pi \frac{lm}{p}} = (-1/p^{n+1})(K'+1)\alpha \ ,
\end{aligned}$$

avec $\alpha = \sum_m \Psi_A(m) \times \begin{cases} -1 \text{ si } m \neq 0 \bmod p \\ p - 1 \text{ si } p|m \end{cases}$. On a:

$$\begin{aligned}
\alpha &= -\sum_m \Psi_A(m) + p \sum_{p|m} \Psi_A(m) \\
\alpha &\leq -(K+1) + p \left( (2/K) \sum_{\lambda=0}^{E(\frac{K}{p})} (K - \lambda p) \right),
\end{aligned}$$

où $E(\frac{K}{p})$ désigne la partie entière de $\frac{K}{p}$ ;

$$\begin{aligned}
\alpha &\leq -(K+1) + (2p/K) \left[ K(E(K/p)+1) - \frac{p}{2} E(K/p)(E(K/p)+1) \right] \\
\alpha &\leq -(K+1) + (2p/K) \left[ (E(K/p)+1)(K - \frac{p}{2} E(K/p)) \right] \\
\alpha &\leq -(K+1) + (2p/K) \left[ (E(K/p)+1)(K - \frac{p}{2}(K/p - 1)) \right] \\
\alpha &\leq -(K+1) + (2p/K)(E(K/p)+1) \left[ (K+p)/2 \right] \\
\alpha &\leq -(K+1) + (K+p)^2/K \\
\alpha &\leq -K - 1 + K + 2p + (p^2/K) \leq p(2 + (p/K)) - 1 \ ,
\end{aligned}$$

or les hypothèses du lemme donnent $K \geq \frac{|A|-1}{2} \geq p$, donc

$$\Lambda_{0,} \geq \frac{-1}{p^{n+1}}(K'+1)(3p-1), \text{ soit } \Lambda_{0,} \geq \frac{-3}{q}(K'+1) .$$

Somme des termes en $(h \neq 0, h' \neq 0)$

$$\Lambda_, = \frac{1}{q^2} \sum_{h=1}^{q-1} \sum_{h'=1}^{q-1} \left( \sum_m \Psi_A(m) e^{\frac{2i\pi hm}{q}} \right) \left( \sum_{m'} \Psi_B(m') e^{\frac{2i\pi h'm'}{q}} \right) S(-h, -h'; q)$$

$$|\Lambda_,| \leq \frac{2\sqrt{2}}{q^2} q^{\frac{1}{2}} \sum_{h=1}^{q-1} \sum_{h'=1}^{q-1} |\sum_m \Psi_A(m) e^{\frac{2i\pi hm}{q}} \sum_{m'} \Psi_B(m') e^{\frac{2i\pi h'm'}{q}}|.(h, h', q)^{\frac{1}{2}}$$

On estime alors, pour $h \neq 0$:

$$\Theta(h) = |\left( \sum_m \Psi_A(m) e^{\frac{2i\pi hm}{q}} \right)|$$

$$\begin{aligned}
\Theta(h) &= \frac{1}{K} |\sum_{k=1}^{K} k e^{2i\pi h \frac{(a+k)}{q}} + \sum_{k=1}^{K} k e^{2i\pi h \frac{(a+2K+1-k)}{q}}| \\
&= \frac{1}{K} |\sum_{k=1}^{K} k \left( e^{2i\pi(-K-\frac{1}{2}+k)\frac{h}{q}} + e^{2i\pi(K+\frac{1}{2}-k)\frac{h}{q}} \right)| \\
&= \frac{2}{K} |\sum_{k=1}^{K} k \cos(2\pi \frac{h}{q}(-K - \frac{1}{2} + k))|
\end{aligned}$$

et par transformée d'Abel,

$$
\begin{aligned}
&= \frac{2}{K}\left|\left(\sum_{k=0}^{K}(k-(k+1))\sum_{l=0}^{k}\cos(\frac{2\pi h}{q}(-K-\frac{1}{2}+l))\right)+(K+1)\left(\sum_{l=0}^{K}\cos(\frac{2\pi h}{q}(-K-\frac{1}{2}+l))\right)\right| \\
&= \frac{2}{K}\left|\sum_{k=0}^{K}\left(\sum_{l=0}^{K}\cos(\frac{2\pi h}{q}(-K-\frac{1}{2}+l))-\sum_{l=0}^{k}\cos(\frac{2\pi h}{q}(-K-\frac{1}{2}+l))\right)\right| \\
&= \frac{2}{K}\left|\sum_{k=0}^{K-1}\left(\sum_{l=k+1}^{K}\cos(\frac{2\pi h}{q}(-K-\frac{1}{2}+l))\right)\right| \\
&= \frac{2}{K}\left|\sum_{k=0}^{K-1}\mathrm{Re}\left(\sum_{l=k+1}^{K}e^{2i\pi\frac{h}{q}(-K-\frac{1}{2}+l)}\right)\right| \\
&= \frac{2}{K}\left|\sum_{k=0}^{K-1}\mathrm{Re}\left(e^{\frac{-i\pi h}{q}}\left(\sum_{l=0}^{K-k-1}e^{-2i\pi\frac{h}{q}l}\right)\right)\right| \\
&= \frac{2}{K}\left|\sum_{k=0}^{K-1}\mathrm{Re}\left(e^{-\frac{i\pi h}{q}}(1-e^{-2i\pi\frac{h}{q}(K-k)})/(1-e^{-2i\pi\frac{h}{q}})\right)\right| \\
&= \frac{2}{K}\left|\sum_{k=0}^{K-1}\mathrm{Re}\left[e^{-i\pi\frac{h}{q}}(e^{-i\pi\frac{h}{q}(K-k)}\cdot\sin(\frac{\pi h}{q}(K-k)))\right]/(e^{-i\pi\frac{h}{q}}\cdot\sin(\frac{\pi h}{q}))\right| ;
\end{aligned}
$$

$$
\begin{aligned}
\Theta(h) &= \frac{2}{K}\left|\sum_{k=0}^{K-1}\left[\sin(\frac{\pi h}{q}(K-k))/\sin(\frac{\pi h}{q})\right]\cos(\frac{\pi h}{q}(K-k))\right| \\
&= 1/(K|\sin\frac{\pi h}{q}|)\cdot\left|\sum_{k=0}^{K-1}\sin\frac{2\pi h}{q}(K-k)\right| \\
&= 1/(K|\sin\frac{\pi h}{q}|)\cdot\left|\mathrm{Im}\left(e^{2i\pi\frac{h}{q}K}\sum_{k=0}^{K-1}e^{-2i\pi\frac{h}{q}k}\right)\right| \\
&= 1/(K|\sin\frac{\pi h}{q}|)\cdot\left|\mathrm{Im}\left(e^{2i\pi\frac{h}{q}K}(1-e^{-2i\pi\frac{h}{q}K})/(1-e^{-2i\pi\frac{h}{q}})\right)\right| \\
&= 1/(K|\sin\frac{\pi h}{q}|)\cdot\left|\mathrm{Im}\left(e^{2i\pi\frac{h}{q}K}(e^{-i\pi\frac{h}{q}K}\sin(\frac{\pi h}{q}K))/(e^{-i\pi\frac{h}{q}}\sin(\frac{\pi h}{q}))\right)\right| ,
\end{aligned}
$$

donc on a enfin :
$$
\Theta(h) = \frac{1}{K}\frac{|\sin(\frac{\pi h}{q}K)\cdot\sin(\frac{\pi h}{q}(K+1))|}{|\sin(\frac{\pi h}{q})|^2} .
$$

On veut majorer $\Theta(h)$ ; on commence donc par $|\frac{\sin(Kx)}{\sin x}| = A_K(x)$, pour $x \in\ ]0,\frac{\pi}{2}]$. On a $A_K(x) \leq K\frac{\pi}{2}|\frac{\sin(Kx)}{Kx}|$ ; or on vérifie que $|\frac{\sin(X)}{X}| \leq \frac{2}{X+1}$ si $X \in\ ]0,\infty]$, donc $A_K(x) \leq \frac{\pi K}{(X+1)}$, soit

$$
A_K(x) \leq \frac{\pi}{x+\frac{1}{K}} \text{ , pour } x \in ]0,\frac{\pi}{2}] .
$$

On en déduit que, pour $h \in \{0, ..., E(\frac{q}{2})\}$, $\Theta(h) \leq \frac{1}{K}\left(\pi/(\frac{\pi h}{q} + \frac{1}{K})\right) \cdot \left(\pi/(\frac{\pi h}{q} + \frac{1}{K+1})\right)$,
soit :

$$\Theta(h) \leq \frac{1}{K}\left(\frac{\pi}{\frac{\pi h}{q} + \frac{1}{K+1}}\right)^2.$$

On revient maintenant à la majoration de $|\Lambda_{,}|$ :

$|\Lambda_{,}| \leq 2\sqrt{2}\, q^{\frac{-3}{2}} \sum_{h=1}^{q-1} \sum_{h'=1}^{q-1} \Theta(h).\Theta(h').(h,h',q)^{\frac{1}{2}}$, or $\Theta(h)$ et $(h,h',q)$ sont invariants par la transformation $h \rightsquigarrow q - h$, donc

$$|\Lambda_{,}| \leq 8\sqrt{2}\, q^{\frac{-3}{2}} \sum_{h=1}^{E(\frac{q}{2})} \sum_{h'=1}^{E(\frac{q}{2})} \frac{\pi^2}{K}\left(\frac{1}{\frac{\pi h}{q} + \frac{1}{K+1}}\right)^2 \cdot \frac{\pi^2}{K'}\left(\frac{1}{\frac{\pi h'}{q} + \frac{1}{K'+1}}\right)^2 (h,h',q)^{\frac{1}{2}},$$

$$|\Lambda_{,}| \leq 8\sqrt{2}\, \pi^4(q^{\frac{-3}{2}}/KK') \sum_{k=0}^{n-1} p^{\frac{k}{2}} \sum_{\substack{h=1 \\ p^k|h}}^{q} \sum_{\substack{h'=1 \\ p^k|h'}}^{q} \left(\frac{1}{\frac{\pi h}{q} + \frac{1}{K+1}}\right)^2 \left(\frac{1}{\frac{\pi h'}{q} + \frac{1}{K'+1}}\right)^2$$

$|\Lambda_{,}| \leq 8\sqrt{2}\, \pi^4(q^{\frac{-3}{2}}/KK') \sum_{k=0}^{n-1} p^{\frac{k}{2}} \sum_{h=1}^{p^{n-k}} \sum_{h'=1}^{p^{n-k}} \left(1/(\frac{\pi h}{p^{n-k}} + \frac{1}{K+1})\right)^2 \cdot \left(1/(\frac{\pi h'}{p^{n-k}} + \frac{1}{K'+1})\right)^2.$

Pour majorer ceci, on considère la fonction $f(x) = \left(1/(\pi x + \frac{1}{K+1})\right)^2$. Comme elle est convexe sur $[0,1]$, on a :

$$\begin{aligned}
\Delta &:= \sum_{h=1}^{p^{n-k}} f(h/p^{n-k}) = p^{n-k} \sum_{h=1}^{p^{n-k}} \frac{1}{p^{n-k}} f(h/p^{n-k}) \\
&\leq p^{n-k} \int_0^1 f(x)dx \leq p^{n-k}\left[\frac{-1}{\pi}\frac{1}{(\pi x + \frac{1}{K+1})}\right]_0^1 \leq p^{n-k}.(K+1)/(\pi + \frac{1}{K+1}).
\end{aligned}$$

D'où :

$$\begin{aligned}
|\Lambda_{,}| &\leq 8\sqrt{2}\, \pi^4(q^{\frac{-3}{2}}/KK') \sum_{k=0}^{n-1} p^{\frac{k}{2}} \left[p^{n-k}(K+1)/(\pi + \frac{1}{K+1})\right] \cdot \left[p^{n-k}(K'+1)/(\pi + \frac{1}{K'+1})\right] \\
|\Lambda_{,}| &\leq 8\sqrt{2}\, \pi^4(p^{\frac{-3n}{2}}/KK') \left[(K+1)(K'+1)/\left((\pi + \frac{1}{K+1})(\pi + \frac{1}{K'+1})\right)\right] \sum_{k=0}^{\infty} p^{2n-\frac{3k}{2}} \\
|\Lambda_{,}| &\leq \left[(8\sqrt{2}\, \pi^2)/\left((1 + \frac{1}{\pi(K+1)})(1 + \frac{1}{\pi(K'+1)})\right)\right] \times (1 + \frac{1}{K})(1 + \frac{1}{K'}).\sqrt{q}/(1 - p^{\frac{-3}{2}}).
\end{aligned}$$

Par hypothèse, $K, K' \geq 5$ (car $\frac{|A|-1}{2} \geq 5$), donc $1 + \frac{1}{K^{(\prime)}} \leq \sqrt{2}$. On en déduit enfin :

$$|\Lambda_{,}| \leq 16\pi^2 \frac{\sqrt{2}}{1 - 2^{\frac{-3}{2}}}\sqrt{q} \leq \frac{64\pi^2}{2\sqrt{2} - 1}\sqrt{q}.$$

On rassemble les estimations de $\Lambda_{0,0}$, $\Lambda_{0,}$ et $\Lambda_{,}$, et on écrit :

$$\begin{aligned}
\Lambda &\geq \frac{1}{2q}(K+1)(K'+1) - \frac{3}{q}\left((K+1) + (K'+1)\right) - \frac{64\pi^2}{2\sqrt{2}-1}\sqrt{q} \\
&\geq \frac{1}{2q}\left[(K-5)(K'-5) - 36\right] - \frac{64\pi^2}{2\sqrt{2}-1}\sqrt{q},
\end{aligned}$$

donc $\Lambda > 0$ pour:

$$(K-5)(K'-5) > 36 + \frac{128\pi^2}{2\sqrt{2}-1}q^{\frac{3}{2}} .$$

Or on a par hypothèse: $(A-11)(B-11) > 144 + (512\pi^2)/(2\sqrt{2}-1)q^{\frac{3}{2}}$, et $K \geq \frac{|A|-1}{2} \geq 5$ (de même pour $K'$), donc

$$(K-5)(K'-5) \geq \frac{1}{4}(|A|-11)(|B|-11) > 36 + \frac{128\pi^2}{2\sqrt{2}-1}q^{\frac{3}{2}} ,$$

ce qu'il fallait démontrer. $\square$

On montre pour finir comment on passe du Lemme 7 à la Proposition 5 : on a vu qu'on pouvait prendre $|A| \geq (p^n/D^2) - 2$ et $B \geq (p^n/D) - D - 2$. Pour satisfaire aux conditions du lemme, il suffit de :

$$\begin{cases} (p^n/D^2) - 2 \geq \max(11, 2p+1), \ (p^n/D) - D - 2 \geq \max(11, 2p+1), \text{ et} \\ ((p^n/D^2) - 13)((p^n/D) - D - 13) > 144 + \frac{512\pi^2}{2\sqrt{2}-1}p^{3n/2} . \end{cases}$$

Si on a les conditions $\begin{cases} D \geq 2 \ ; \ p^n/D^2 \geq \sup(26, 4p) \ ; \ \text{et} \\ \frac{1}{2}(p^n/D^2).\frac{1}{2}(p^n/D) > 144 + \frac{512\pi^2}{2\sqrt{2}-1}p^{\frac{3n}{2}} , \end{cases}$
alors les précédentes sont satisfaites, car : $p^n \geq 26D^2$ implique $13+D \leq \frac{1}{2}(p^n/D)$ puisque $D \geq 2 \Rightarrow 26D \geq 2D + 26$ ; de plus $p^n/D^2 \geq 4p \Rightarrow (p^n/D^2) - 2 \geq 2p+1$ et $p^n/D \geq 4pD \geq D + 2p + 3$ (en effet, $\frac{D+3}{4D-2} < 2 \leq p$ si $D \geq 2$).

Comme enfin on a $144 \leq \frac{1}{2}\left(\frac{512}{2\sqrt{2}-1}\pi^2 p^{3n/2}\right)$ puisque $p \geq 2$, il suffit pour que les conditions du lemme soient vérifiées qu'on ait :

$$\begin{cases} D \geq 2, \ (p^{n-1}/D^2) \geq 4, \ (p^n/D^2) \geq 26, \text{ et} \\ p^{2n}/D^3 > \frac{4096\pi^2}{2\sqrt{2}-1}p^{3n/2} \end{cases}$$

Mais on ne s'intéresse ici qu'aux puissances supérieures ou égales à 2 des nombres premiers ; comme $C = \frac{4096\pi^2}{2\sqrt{2}-1}$ est de beaucoup supérieure à $\sqrt{26}$ et à 4, la dernière condition implique les deux prédentes dans la dernière accolade, et les conditions du lemme sont remplies pour

$$q > \left(\frac{4096\pi^2}{2\sqrt{2}-1}\right)^2 D^6 \ \square.$$

# Références

[Edixhoven 94]  B. Edixhoven, *Rationnal torsion points on elliptic curves over number fields*, Astérisque (Séminaire Bourbaki 782) (1993).

[Hooley]  Hooley, *Applications of Sieve methods to Number Theory*, Cambridge University Press, 1935, page 35.

[Kamienny 92a]  S. Kamienny, *Torsion points on elliptic curves and q-coefficients of modular forms*, Inventiones Mathematicae 109 (1992), pages 221-229

[Kamienny 92b]  S. Kamienny, *Torsion points on elliptic curves over fields of higher degree*, International Mathematics Research Notices 6 (1992).

[Kamienny-Mazur]  S. Kamienny et B. Mazur, *Rationnal torsion of prime order in elliptic curves over number fields*, Astérisque (à paraître)

[Manin 69]  Manin, *A Uniform Bound for p-torsion of Elliptic Curves*, Izv. Akad. Nau. CCCP 33 (1969)

[Manin 72]  Y. Manin, *Parabolic points and zeta function of modular curves*, Math. USSR Izvestija 6 (1972), pages 19-64

[Mazur 77]  B. Mazur, *Modular curves and the Eisenstein ideal*, Publications mathématiques de l'I.H.E.S. 47 (1977), pages 33-186 .

[Merel 93]  L. Merel, *Sur quelques aspects géométriques et arithmétiques de la théorie des symboles modulaires*, Thèse de doctorat de l'université de Paris VII (1993)

[Merel 95]  L. Merel, *Bornes pour la torsion des courbes elliptiques sur les corps de nombres*, Inventiones Mathematicae, à paraître.

[Oesterlé 94]  J. Oesterlé, *Torsion des courbes elliptiques sur les corps de nombres*, article?? à paraître

[Salié]  Salié, *Über die Kloostermanschen Summen*, Math. Zeit. 34, 1931, pages 91-109.

---

Pierre Parent

ENS

45, rue d'Ulm

75005 Paris

E-mail: `parent@clipper.ens.fr`